\documentclass[10pt,twoside]{article}
\usepackage{amssymb, amsmath, latexsym, a4}
\usepackage[all]{xy}
\def\YEAR{\year}\newcount\VOL\VOL=\YEAR\advance\VOL by-1995
\def\firstpage{1}\def\lastpage{1000}
\def\received{}\def\revised{}
\def\communicated{}
\makeatletter
\def\magnification{\afterassignment\m@g\count@}
\def\m@g{\mag=\count@\hsize6.5truein\vsize8.9truein\dimen\footins8truein}
\makeatother
\font\eightrm=cmr8
\font\caps=cmcsc10                    
\font\Caps=cmcsc10 scaled \magstep1   

\makeatletter
\@specialpagefalse\headheight=8.5pt
\def\DocMath{}
\renewcommand{\@evenhead}{%
    \ifnum\thepage>\lastpage\rlap{\thepage}\hfill%
    \else\rlap{\thepage}\slshape\leftmark\hfill{\caps\SAuthor}\hfill\fi}%
\renewcommand{\@oddhead}{%
    \ifnum\thepage=\firstpage{\DocMath\hfill\llap{\thepage}}%
    \else{\slshape\rightmark}\hfill{\caps\STitle}\hfill\llap{\thepage}\fi}%
\makeatother
\def\TSkip{\bigskip}
\newbox\TheTitle{\obeylines\gdef\GetTitle #1
\ShortTitle  #2
\SubTitle    #3
\Author      #4
\ShortAuthor #5
\EndTitle
{\setbox\TheTitle=\vbox{\baselineskip=20pt\let\par=\cr\obeylines%
\halign{\centerline{\Caps##}\cr\noalign{\medskip}\cr#1\cr}}%
    \copy\TheTitle\TSkip\TSkip%
\def\next{#2}\ifx\next\empty\gdef\STitle{#1}\else\gdef\STitle{#2}\fi%
\def\next{#3}\ifx\next\empty%
    \else\setbox\TheTitle=\vbox{\baselineskip=20pt\let\par=\cr\obeylines%
    \halign{\centerline{\caps##} #3\cr}}\copy\TheTitle\TSkip\TSkip\fi%
\centerline{\caps #4}\TSkip\TSkip%
\def\next{#5}\ifx\next\empty\gdef\SAuthor{#4}\else\gdef\SAuthor{#5}\fi%
\ifx\received\empty\relax
    \else\centerline{\eightrm Received: \received}\fi%
\ifx\revised\empty\TSkip%
    \else\centerline{\eightrm Revised: \revised}\TSkip\fi%
\ifx\communicated\empty\relax
    \else\centerline{\eightrm Communicated by \communicated}\fi\TSkip\TSkip%
\catcode'015=5}}\def\Title{\obeylines\GetTitle}
\def\Abstract{\begingroup\narrower
    \parskip=\medskipamount\parindent=0pt{\caps Abstract. }}
\def\EndAbstract{\par\endgroup\TSkip}
\long\def\MSC#1\EndMSC{\def\arg{#1}\ifx\arg\empty\relax\else
     {\par\narrower\noindent%
     2000 Mathematics Subject Classification: #1\par}\fi}
\long\def\KEY#1\EndKEY{\def\arg{#1}\ifx\arg\empty\relax\else
    {\par\narrower\noindent Keywords and Phrases: #1\par}\fi\TSkip}
\newbox\TheAdd\def\Addresses{\vfill\copy\TheAdd\vfill
    \ifodd\number\lastpage\vfill\eject\phantom{.}\vfill\eject\fi}
{\obeylines\gdef\GetAddress #1
\Address #2
\Address #3
\Address #4
\EndAddress
{\def\xs{4.3truecm}\parindent=0pt
\setbox0=\vtop{{\obeylines\hsize=\xs#1\par}}\def\next{#2}
\ifx\next\empty 
     \setbox\TheAdd=\hbox to\hsize{\hfill\copy0\hfill}
\else\setbox1=\vtop{{\obeylines\hsize=\xs#2\par}}\def\next{#3}
\ifx\next\empty 
     \setbox\TheAdd=\hbox to\hsize{\hfill\copy0\hfill\copy1\hfill}
\else\setbox2=\vtop{{\obeylines\hsize=\xs#3\par}}\def\next{#4}
\ifx\next\empty\ 
     \setbox\TheAdd=\vtop{\hbox to\hsize{\hfill\copy0\hfill\copy1\hfill}
                \vskip20pt\hbox to\hsize{\hfill\copy2\hfill}}
\else\setbox3=\vtop{{\obeylines\hsize=\xs#4\par}}
     \setbox\TheAdd=\vtop{\hbox to\hsize{\hfill\copy0\hfill\copy1\hfill}
            \vskip20pt\hbox to\hsize{\hfill\copy2\hfill\copy3\hfill}}
\fi\fi\fi\catcode'015=5}}\gdef\Address{\obeylines\GetAddress}
\begin{document}
\Title Comparison of abelian categories recollements
\ShortTitle Recollements comparison
\SubTitle

\Author Vincent FRANJOU\footnote{membre du laboratoire Jean-Leray, UMR 6629 UN/CNRS
} and Teimuraz PIRASHVILI\footnote{supported by the grants INTAS-99-00817 and RTN-Network ``K-theory, linear algebraic groups and related structures'' HPRN-CT-2002- 00287}
\ShortAuthor V. Franjou and T. Pirashvili
\EndTitle
\Abstract
We give a necessary and sufficient condition for a morphism between
recollements of abelian categories to be an equivalence.


\EndAbstract
\MSC 18F, 18E40, 16G, 16D90, 16E
\EndMSC
\KEY
recollement, abelian category, functor
\EndKEY
\Address Universit\'e de Nantes\\d\'ep. de math\'ematiques
2, rue de la Houssini\`ere\\BP 92208\\44322 Nantes cedex 3\\France
\Address A. M. Razmadze\\Mathematical Institute
Aleksidze str.1\\
Tbilisi 380093\\Republic of Georgia
\Address
\Address
\EndAddress
\renewcommand{\theenumi}{\roman{enumi}}
\newtheorem{De}{Definition}[section]
\newtheorem{Th}[De]{Theorem}
\newtheorem{Prop}[De]{Proposition}
\newtheorem{Le}[De]{Lemma}
\newtheorem{Co}[De]{Corollary}
\newcommand{\rdg}{\hfill $\Box $\par\medskip}
\def\t{\otimes }
\def \x{\times}
\def \Ker{{\rm Ker}\ }
\def \Coker{{\rm Coker}\ }
\def \Im{{\rm Im}\ }
\def \Hom{{\rm Hom}}
\def \Ext{{\rm Ext}}
\def \h{{\rm H}}
\def \L{{\rm L}}
\def \R{{\rm R}}
\newcommand{\FF}{\mathbb{F}}
\def \A{{\cal A}}
\def \B{{\cal B}}
\def \C{{\cal C}}
\def \F{{\cal F}}
\def \G{{\cal G}}
\def \Ab{{{\cal A}b}}
\def \Su{{\cal M}}
\def \qh{\text{pre-hereditary}}
\section{Introduction}
Recollements of abelian and triangulated categories
play an important role in geometry of singular spaces
\cite{bbd}, in representation theory
\cite{cps, ps}, in polynomial functors theory
\cite{kuhn, kuhnstrat,rmi} and in ring theory,
where recollements are known as torsion, torsion-free theories \cite{J}.
A fundamental example of recollement of abelian categories
is due to MacPherson and Vilonen \cite{MV}.
It first appeared as an inductive step in the construction of perverse sheaves.
The main motivation for our work was to understand when
a recollement can be obtained through the construction of MacPherson
and Vilonen.

\par\medskip
A {\it recollement situation} consists of three abelian categories
${\A }'$, ${\A }$, ${\A }''$ together with additive functors:
$$\begin{array}{rcccl}
 \ &\ \buildrel{i^*}\over{\longleftarrow}&&\buildrel {j_!}\over\longleftarrow&\\
\A' &\buildrel {i_*}\over\longrightarrow& {\A}&\buildrel {j^*}\over\longrightarrow&{\A}''\\
\ & \buildrel {i^!}\over \longleftarrow& &\buildrel
{j_*}\over\longleftarrow &
\end{array}$$
which satisfy the following conditions:
\begin{enumerate}
    \item  $j_!$ is left adjoint to $j^*$ and $j^*$ is left adjoint to $j_*$
    \item  the unit $Id_{{\A }''}\to j^*j_!$ and the counit  $j^*j_*\to
  Id_{{\A }''}$ are isomorphisms
    \item  $i^*$ is left adjoint to $i_*$ and $i_*$ is left adjoint to $i^!$
    \item  the unit $Id_{{\A }'}\to i^!i_*$ and the counit
  $i^*i_*\to Id_{{\A }'}$ are isomorphisms
    \item  $i_*$ is an embedding onto the full subcategory of $\A $ with
objects $A$ such that $j^*A=0$.
\end{enumerate}
In this case one says that $\A $ is a \emph{recollement of
${\A }''$ and ${\A }'$.} \emph{These notations will be kept throughout the paper.}
Thus in any recollement situation, the category $i_*\A'$ is a localizing and colocalizing
subcategory of $\A$ in the sense of \cite{gab}, and
the category $\A''$ is equivalent to the quotient category of $\A$ by $i_*\A'$.
\par\medskip
If $\B $ is also a recollement of ${\A }''$ and ${\A }'$, then a
\emph{comparison functor} $\A\to \B$ is an exact functor which
commutes with all the structural functors $i^*, i_*,i^!,j_!,j^*,j_*$.
According to \cite[Theorem 2.5]{ps}, a comparison functor between recollements of
triangulated categories is an equivalence of categories. Our example in
Section \ref{contrex} shows that this is not necessarily the case for recollements of abelian categories.

Our main result, Theorem \ref{comparison}, characterizes which
comparisons of recollements are equivalences of categories. As an
application, we give a homological criterion deciding when a
recollement can be obtained through the construction of MacPherson
and Vilonen.

\textsc{Theorem.}\label{mainMV} \emph{A recollement situation of categories with enough projectives
 is isomorphic to a MacPherson-Vilonen
construction if and only if the following two conditions hold.
\begin{enumerate}
\item There exists an exact functor $r$: $\A \to \A'$ such that
$r\circ i_*=Id_{\A '}$ .
\item
For any projective object $V$ of the category $\A'$,
 $(\L _2i^*)(i_*V)=0$.
\end{enumerate}}
\section{Examples}\label{ex}
Our examples are related to polynomial functors.
The relevance of this formalism to polynomial functors was
stressed by N. Kuhn \cite{kuhn}.
\par
We let $\A'$ be the category of finite
vector spaces over the field with two elements ${\FF}_2$, and we
let $\A''$ be the category of finite vector spaces over ${\FF}_2$ with involution, or finite representations of $\Sigma_2$
over ${\FF}_2$.
\subsection{}\label{1stex}
In the first example, the category $\A$ is a category of diagrams of
finite vector spaces over ${\FF}_2$: $$(V_1,H,V_2,P):\ V_1
\rightleftarrows V_2\ ,$$ where $H$: $V_1\to V_2$ and $P$: $V_2\to
V_1$ are linear maps which satisfy: $PHP=0$ and $HPH=0$. The
category $\A$ is equivalent to the category of quadratic functors
from finitely generated free abelian groups to vector spaces over
${\FF}_2$.
 It is a recollement for the following functors:
 $$i^*(V_1,H,V_2,P)={\Coker}(P),\ j_!(V,T)=(V_T,1+T,V,p)$$
 $$i_*(V)=(V,0,0,0), \ \ \ \ \ j^*(V_1,H,V_2,P)= (V_2,HP-1)$$
$$i^!(V_1,H,V_2,P)=\Ker(H)\ ,\ \ j_*(V,T)=(V^T,h,V,1+T)\ ,$$
where $V^T=\Ker (1-T)$, $V_T=\Coker (1-T)$, $h$ is the inclusion
and $p$ is the quotient map. Note that the functor $i_*$ admits an
obvious exact retraction $r$: $(V_1,H,V_2,P)\mapsto V_1$.
\subsection{Comparison fails for abelian categories recollements}\label{contrex}
We now consider the full subcategory of the category $\A$ in Example \ref{1stex},
whose objects satisfy the relation: $PH=0$. This
category is equivalent to the category of quadratic functors from
finite vector spaces to vector spaces over ${\FF}_2$. The same
formulae define a recollement as well. As a result, the inclusion
of categories is a comparison functor. It is not, however, an
equivalence of categories.
\par
\section{The construction of MacPherson and Vilonen \cite{MV}}  \label{macp}
\subsection{} Let $\A'$ and
$\A''$ be abelian categories. Let $F$: $\A''\to \A'$ be a right
exact functor, let $G$: $\A''\to \A'$ be a left exact functor and
let $\xi$: $F\to G$ be a natural transformation. Define the
category $\A(\xi)$ as follows. The objects of $\A(\xi)$ are tuples
$(X,V,\alpha,\beta)$, where $X$ is in $\A''$, $V$ is in $\A'$,
$\alpha:F(X)\to V$ and $\beta:V\to G(X)$ are morphisms in $\A'$ such
that the following diagram commutes:
$$\xymatrix{
F(X)\ar[rr]^{\xi_X}\ar[rd]_{\alpha}&& G(X)\\
&V\ar[ru]_{\beta}& }\ .$$ A morphism from $(X,V,\alpha,\beta)$ to
$(X',V',\alpha',\beta')$ is a pair $(f,\varphi)$, where $f:X\to X'$ is
a morphism in $\A''$ and $\varphi$: $V\to V'$ is a morphism in
$\A'$, such that the following diagram commutes:
$$\xymatrix{
F(X)\ar[r]^{\alpha}\ar[d]^{F(f)}& V\ar[r]^{\beta}\ar[d]^{\varphi}&
G(X)\ar[d]^{G(f)}\\
F(X')\ar[r]^{\alpha'}& V'\ar[r]^{\beta'}& G(X')\\ }\ .$$ The
category $\A(\xi)$ comes with functors:
$$i^*(X,V,\alpha,\beta)={\Coker}(\alpha)\ , \qquad j_!(X)=(X,F(X),Id_{F(X)},\xi_X)\ ,$$
$$i_*(V)=(0,V,0,0)\ , \qquad j^*(X,V, \alpha, \beta)=X\ ,$$
$$i^!(X,V,\alpha,\beta)=\Ker(\beta)\ ,\qquad j_*(X)=(X,G(X),\xi_X,Id_{G(X)}) \ .$$
The functor $i_*$ has a retraction functor $r$:
$$r(X,V,\alpha,\beta)=V\ .$$
The category $\A(\xi)$ is abelian in such a way that the functors
$r$ and $j^*$ are exact. The above data define a recollement. Note
that we recover the natural transformation $\xi$ from the
retraction $r$ and the recollement data as:
$$F=rj_!\qquad G=rj_*\qquad  \xi\simeq rN\ .$$
The category $\A$ depends only \cite[Proposition 1.2]{MV}
on the class of the extension
$$0\to i^!j_!\to  F\buildrel \xi \over\rightarrow
G\to i^*j_*\to 0\ ,$$ image by $r$ of the exact sequence
(\ref{norme}).
\subsection{}\label{semidirect}
We now consider two particular cases of this construction,
already known to Grothendieck (see \cite{A}).
Let $F$: $\A''\to \A'$ be a right exact functor. Take $\xi$: $F\to 0$ to be
the transformation into the trivial functor. The corresponding
construction is denoted by ${\A '}\rtimes _F{\A ''}$. Thus objects
of this category are triples $(V,X,\alpha)$, where $V$ and $X$ are
objects of $\A '$ and $\A ''$ respectively and $\alpha$ is a morphism
$\alpha$: $F(X)\to V$ of the category $\A '$. Note that
$i^*j_*=0$ and $i^!j_!\cong F$. Moreover, $i^!$ and $j_*$ are
exact functors.
\par
Similarly, let  $\B'$ and $\B''$ be  abelian categories and let
$G$: $\B''\to \B'$ be a left exact functor.
We take $\xi$: $0\to G$ to be the natural
transformation from the trivial functor. The corresponding
recollement is denoted by ${\cal B'}\ltimes _G{\cal B''}$. Objects of
this category are triples $(B'',B',\beta: B'\to G(B''))$.
Assuming now $\B'=\A''$, $\B''=\A'$ and $G:\A'\to \A''$ is right
adjoint to $F$, the category ${\A '}\rtimes _F{\A ''}={\cal
A''}\ltimes _G\A'$ fits into two different recollement situations.
\section{General properties of recollements}
 Most of the properties in this section can probably be found in
\cite{bbd} or other references. We list them for convenience.
Note however that, when they are not a consequence of
\cite{gab}, they are usually stated and
proved in the context of triangulated categories. We consistently
provide statements (and a few proofs) in the context of abelian
categories and derived functors.
\subsection{First properties}\label{properties}
We remark as usual that taking opposite categories results in the
exchange of $j_!$ and $ i^*$ with $j_*$ and $i^!$ respectively.
This is referred to as duality. For instance, the relation
$j^*i_*=0$ - a consequence of (v) - yields the dual relation
$i^!j_*=0$.
\begin{Prop}\label{nulia} In any recollement situation:
$$i^*j_!=0\ , \ \ i^!j_*=0.$$
\end{Prop}
\begin{Prop}
The units and counits of adjonction
give rise to exact sequences of natural transformations:
\begin{equation}\label{epsilon}
  j_!j^* \buildrel \epsilon \over\rightarrow Id_{\A }\to i_*i^*\to 0
 \end{equation}
\begin{equation}\label{eta}
 0\to i_*i^!\to Id_{\A } \buildrel \eta \over\rightarrow j_*j^* \ .
 \end{equation}
\end{Prop}
\medskip
We now recall the definition of the norm $N$: $j_!\to j_*$. For
any $X$, $Y$ in $\A ''$, there are natural isomorphisms:
$$\Hom_{\A }(j_!X,j_*Y)\cong \Hom_{\A ''}(X,j^*j_*Y)\cong \Hom_{\cal
A''}(X,Y).$$ For $Y=X$, let $N_X$: $j_!X\to j_*X$ be the map
corresponding to the identity of $X$.
It is a natural transformation
\cite[1.4.6.2]{bbd}. The norm $N$ is thus defined so that:
$Nj^*=\eta\circ\epsilon$. Hence:
\begin{equation}\label{norm}
N\cong N(j^*j_*)=(Nj^*)j_*\cong(\eta\circ\epsilon)j_*= \eta
j_*\circ\epsilon j_*\cong\epsilon j_*\ \text{and, dually}\ N\cong\eta j_!\ .\
\end{equation}
The image of the norm is a functor
 $$j_{!*}:=\Im(N:j_!\to j_*):\ \A''\to\A\ .$$
\begin{Prop}\label{interm} In any recollement situation:
 $i^*j_{!*}=0$ , $\ i^!j_{!*}=0$ .
\end{Prop}
{\it Proof}. Use Proposition \ref{nulia} and apply $i^*$ to the epi $j_!\to j_{!*}$.\rdg
\begin{Prop}\label{birtvi}
In any recollement situation, there is a short
exact sequence of natural transformations
\begin{equation} \label{norme}
0\to i_*i^!j_!\to j_!\buildrel N\over\rightarrow j_*\to
i_*i^*j_*\to 0\ .
\end{equation}
\end{Prop}
{\it Proof}. Precompose the exact sequence
(\ref{epsilon}) with $j_*$. Precomposition is exact, hence one
gets the following exact sequence:
$$j_! \to j_* \to i_*i^*j_*\to 0\ ,$$
where the left arrow is the norm $N$ according to (\ref{norm}).
Dually, there is an exact sequence:
$$ 0\to i_*i^!j_!\to j_!\buildrel N\over\rightarrow j_*\ .$$
Splicing the two sequences together gives the result. \rdg
\par
Applying the snake lemma, 
one gets the following strong restriction on the functors $i^!j_!$
and $i^*j_*$ of a recollement situation.
\begin{Co}
For any short exact sequence in $\A''$:
$$0\to X\to Y \to Z\to 0$$
there is an exact sequence in $\A'$:
$$i^!j_!(X)\to i^!j_!(Y)\to i^!j_!(Z)\to i^*j_*(X)\to   i^*j_*(Y)\to
 i^*j_*(Z)$$
\end{Co}
\subsection{Homological properties}
In this section we investigate the derived functors of the
functors in a recollement situation. We use the following
convention throughout this section: When mentioning left derived
functors $\L-$, the category $\A$, and thus the categories $\A'$
and $\A''$, have enough projectives, and, similarly, when
mentioning right derived functors $\R-$, the categories $\A$,
$\A'$ and $\A''$ have enough injectives. Most of the proofs
consist in applying long exact sequences for derived functors to
Section \ref{properties}'s exact sequences.
\begin{Prop}\label{stan} For each integer $n\geq 1$:
$$j^*(\L_nj_!)=0\ \ , \ \ \ j^*(\R ^nj_*)=0\ .$$
\end{Prop}
\begin{Prop}
\begin{equation}\label{88}
(\L_1i^*)i_*=0\ \  ,\ \ \ \ (\R ^1i^!)i_*=0
\end{equation}
\begin{equation} \label{pirvelicarmoebuli}
(\L_1i^*)j_!=0\ \  ,\ \ \ \ (\R^1i^!)j_*=0
\end{equation}
\begin{equation} \label{1142}
(\L_1i^*)j_{!*}=i^!j_!\ \  ,\ \ \ \ (\R^1i^!)j_{!*}=i^*j_*
\end{equation}
\end{Prop}
\begin{Prop}\label{spectral} There is a natural exact sequence:
\begin{eqnarray}
\nonumber 0\to \Ext^1_{\A '}(i^*A,V)\to \Ext^1_{\A}(A,i_*V)
\buildrel \eta \over \longrightarrow\Hom_{\A'}((\L_1i^*)A,V)\to \ \ \ \ \\
\nonumber \to\Ext^2_{\A'}(i^*A,V)\to\Ext^2_{\A}(A,i_*V)\ .
\end{eqnarray}
\end{Prop}
{\it Proof}. This follows from the spectral sequence for the derived
functors of the composite functors:
\begin{equation}\label{sps}
E^2_{pq}=\Ext_{\A '}^p(\L _qi^*(A),V)\Longrightarrow
\Ext^{p+q}_{\A }(A,i_*V) \ .
\end{equation}\rdg
\begin{Prop}\label{epsilononKer}
Let $A$ be an object in $\Ker i^*$. The counit $\epsilon_A$:
$j_!j^*A\to A$ is epi and its kernel is in $i_*\A'$. Indeed, if
$\A$ has enough projectives, there is a short exact sequence:
\begin{equation}\label{kerepsilon}
0\to i_*(\L_1i^*)A\to j_!j^*A  \buildrel \epsilon_A
\over\longrightarrow A \to 0\  .
\end{equation}
\end{Prop}
We prove the dual statement:
\begin{Prop}\label{etaonKer}
Let $A$ be an object in $\Ker i^!$. The unit $\eta_A$: $A\to
j_*j^*A$ is mono and its cokernel is in $i_*\A'$. Indeed, if $\A$
has enough injectives, there is a short exact sequence:
\begin{equation}\label{cokereta}
0\to A  \buildrel \eta_A \over\longrightarrow j_*j^*A \to
i_*(\R^1i^!)A\to 0\  .
\end{equation}
\end{Prop}
{\it Proof}. When $i^!A=0$, the exact sequence (\ref{eta})
simplifies to a short exact sequence:
\begin{equation}\label{cokeretaonly}
0\to A\buildrel \eta_A \over\longrightarrow j_*j^*A\to
\Coker\eta_A\to 0\ .
\end{equation}
First applying the exact functor $j^*$, and using that $j^*\eta$
is an iso, we see that $j^*(\Coker\eta_A)=0$. Thus $\Coker\eta_A$
is in $i_*\A '$. Suppose that $\A$ has enough injectives. Applying
now the left exact functor $i^!$, the long exact sequence for
right derived functors gives an exact sequence:
$$0\to i^!A\to i^!j_*j^*A\to i^!\Coker\eta_A\to (\R^1i^!)A\to
(\R^1i^!)j_*j^*A\ .$$ Proposition \ref{nulia} and
(\ref{pirvelicarmoebuli}) give an isomorphism
$i^!\Coker(\eta_A)\cong \R ^1i^!(A)$.\rdg
\subsection{Description of the image of $j_*$, $j_{!*}$, $j_!$ }
Since $j^*j_!\cong j^*j_*\cong j^* j_{!*}\cong Id_{\A''}$,
the functors $j_!,j_*, j_{!*}$: $\A ''\to \A$ are full
embeddings. The next result describes the essential image of each
of them.
\begin{Prop}\label{anasaxebi} The functors $j_!,j_*,j_{!*}:{\A ''}\to \cal
A$ induce the following equivalences of categories:
$$ j_{!*}:{\A }''\to \{A\in {\A }\mid i^*(A)=0=i^!(A)\},$$
$$j_!:{\A }''\to \{A\in {\A }\mid i^*(A)=0=\L _1i^*(A)\},$$
$$j_*:{\A }''\to \{A\in {\A }\mid i^!(A)=0=\R ^1i^!(A)\}.$$
\end{Prop}
\subsection{A monomorphism on $\Ext$-groups}
 Since $j^*:\A \to A''$ is an exact functor, it induces an homomorphism
$$\Ext^n_{\A }(A,B)\to \Ext_{\A ''}^n(j^*A,j^*B), \ n\geq 0.$$
It is well-known that when $A$ and $B$ are simple objects, this map is injective
for $n=1$ (see for example \cite[Proposition 4.12 ]{kuhn}).
The following more general result holds.
\begin{Prop} Let $A,B\in\A $ be objects for which $i^*A=0$ and $i^!B=0$.
Suppose $j^*A\not =0$ and $j^*B\not =0$. Then
$$\Ext^1_{\A }(A,B)\to \Ext_{\A ''}^1(j^*A,j^*B)$$
is a monomorphism.
\end{Prop}
\section{Description of $\Ker i^*$ and $\Ker i^!$}\label{717}
Let $\Ker i^!$ be the full subcategory of objects $A$ of $\A $
such that $i^!A=0$, and let $\Ker i^*$ be the full subcategory of
objects $A$ of $\A $ such that $i^*A=0$. In this section, we
describe these subcategories of $\A$ in terms of the categories
$\A'$, $\A''$, and the functors $i^*j_*$, $i^!j_!$ between them,
through the following construction:
\begin{De}
Let $T$: $\A''\to \A'$ be an additive functor between abelian
categories. The category ${\Su}(T)$ has objects triples
$(X,V,\alpha)$ where $X$ is in $\A''$, $V$ is in $\A'$, and
$\alpha$: $V\to TX$ is a monomorphism. A map from $(X,V,\alpha)$
to $(X',V',\alpha')$ is a pair of morphisms $(f,\varphi)$ such
that the following diagram commutes:
$$\xymatrix{
V\ar[r]^{\alpha}\ar[d]^{\varphi}&
T(X)\ar[d]^{T(f)}\\
V'\ar[r]^{\alpha'}& T(X')\ .\\ }$$
\end{De}
\par\medskip
The following theorem is inspired by \cite{moambe}.
\begin{Th}\label{qartuli}
In a recollement with enough projectives, the functor
$A\mapsto (j^*A,\ i^*A,\ i^*\eta _A: i^*A\to i^*j_*j^*A)$ is an
equivalence from the category $\Ker i^!$ to the category
${\Su}(i^*j_*)$.
\end{Th}
{\it Proof}. First, we show that the functor is well defined.
Apply the functor $i^*$ on the short exact sequence
(\ref{cokeretaonly}). There results an
exact sequence:
$$\L_1i^*(\Coker\eta_A)\to i^*A\to i^*j_*j^*(A)\to i^*\Coker\eta_A\to 0\ .$$
whose left term cancels by Proposition \ref{etaonKer} and (\ref{88}).
The map $i^*\eta _A$ is thus mono.
\par\medskip
Next, we define the quasi-inverse: ${\Su}(i^*j_*)\to\Ker i^!$. To
an object $(X,V,\alpha)$, it associates the kernel $A(X,V,\alpha)$
of the composite of epis:
$$j_*X  \buildrel \epsilon j_* \over\rightarrow i_*i^*j_*X \to
\Coker i_*\alpha\ .$$ That is, $A(X,V,\alpha)$ fits in the
following map of extensions:
$$\xymatrix{0\ar[r]& j_{!*}X \ar[r]& j_*X \ar[r] & i_*i^*j_*X \ar[r]&0\ \ \\
0\ar[r]& j_{!*}X \ar[r]\ar[u]^{=}& A(X,V,\alpha)\ar[r]\ar[u]&i_*V
\ar[r]\ar[u]^{i_*\alpha}&0 \ .}
$$
To a map $(f,\varphi)$, it associates the map induced by $j_*(f)$.
\par\medskip
We leave the verifications to the reader, with the help of the
isomorphism $Nj^*\cong\epsilon\circ\eta$.\rdg
\par
The dual study of the category $\Ker i^*$ leads to the following.
\begin{Th}
In a recollement with enough injectives, the functor
$A\mapsto (j^*A,\ i^!\Ker\epsilon_A,\ i^!\Ker\epsilon_A\to i^!j_!j^*A)$ is an
equivalence from the category $\Ker i^*$ to the category
${\Su}(i^!j_!)$.
\end{Th}
This time, the quasi-inverse fits in the following map of
extensions:
$$\xymatrix{
0\ar[r]& i_*i^!j_!X \ar[r]\ar[d]& j_!X \ar[r]\ar[d] & j_{!*}X \ar[r]\ar[d]^{=}&0\ \ \\
0\ar[r]& \Coker (i_*\alpha)\ar[r]& A(X,V,\alpha) \ar[r] & j_{!*}X
\ar[r]&0\ .}
$$
Note (Proposition \ref{epsilononKer}) that when the recollement has enough projectives,
$i^!\Ker\epsilon_A$ is nothing but $(\L_1i^*)A$.
\section{Recollements as linear extensions}\label{linext}
The exact sequence (\ref{eta}) tells that every object $A$ in $\A$
sits in a short exact sequence:
$$
0\to \Ker\eta_A \to A \buildrel \eta_A \over\rightarrow \Im
\eta_A\to 0\ .
$$
where $\Ker\eta_A\cong i_*i^!A$ is in $i_*\A'$ and
$\Im\eta_A\cong A/i_*i^!A$ is in $\Ker i^!$. We denote by $\G$
the category encoding these data from the recollement
situation. That is, objects of the category $\G$ are triples
$(A,U,e)$ of an object $A$ in $\Ker i^!$, an object $U$ in $\A'$
and an extension class $e$ in the group $\Ext ^1_{\A}(A, i_*U)$. A
map from $(A,U,e)$ to $(A',U',e')$ is a pair of morphism
$(\alpha:A\to A',
\beta:U\to U')$ such that: $\alpha^*e'=(i_*\beta)_*e$ in the group $\Ext
^1_{\A}(A', i_*U)$. It comes with a functor:
$$\A\to\G \ \ \ \ \ \ \ B\mapsto (\Im\eta_B,\ i^!B,\
[0\to \Ker\eta_B \to B \buildrel \eta \over\longrightarrow \Im
\eta_B\to 0])\ .
$$
Because of the Yoneda correspondance between extensions and
elements in $\Ext^1$, this functor induces an equivalence of
categories to $\G$ from the following category $\B$. The objects
of $\B$ are those of $\A$, and a map in $\Hom _{\B}(B,B')$ is a
class of maps in $\Hom_{\A}(B,B')$ inducing the same map in $\G$.
\par\medskip
We claim that $\A\to \B$ defines a linear extension of categories
in the sense of Baues and Wirsching. For completeness, we now
recall what we need from this theory (however, the following
defining properties might be better understood by just looking at
our example).
\begin{De} \rm{\cite[IV.3]{BW}}
Let $\B $ be a category and let $D:{\B } ^{op} \times {\B }
\rightarrow \Ab$ be a bifunctor with abelian groups values. We say
that
\begin{equation}\label{311}
\xymatrix{0\ar[r]& D  \ar[r] & {\C } \ar[r]^{p} &{\B }\ar[r] &0}
\end{equation}
is a  linear extension of the category $\B $ by $D$ if the
following conditions hold:
\begin{enumerate}
\item  $\C $ is a category and $p$ is a functor. Moreover $\C $
and $\B $ have the same objects, $p$ is the identity on objects
and $p$ is surjective on morphisms. \item  For any objects $c$ and
$d$ in $\B$, the abelian group $D (c,d)$ acts on the set $\Hom_{\C
}(c,d)$. Moreover
 $p(f_0)=p(g_0)$ if and only if there is unique $\alpha$ in $D(c,d)$ such
that: $g_0=f_0 +\alpha$. Here for each $f_0:c\rightarrow d$ in $\C
$ and $\alpha \in D(c,d)$ we write $f_0 +\alpha$ for the action of
$\alpha $ on $f_0 $. \item  The action satisfies the linear
distributivity law: for two composable maps $f_0$ and $g_0$ in
$\C$
$$(f_0 +\alpha)(g_0 +\beta) =f_0g_0 +f_* \beta +g^* \alpha \ ,$$
where $f=p(f_0)$ and $g=p(g_0)$.
\end{enumerate}
\end{De}
A morphism between two linear extensions
$$\xymatrix{
0\ar[r]& D  \ar[r]\ar[d]^{\phi_1} & {\C }
\ar[r]^{p}\ar[d]^{\phi_0} &
{\B }\ar[r]\ar[d]^\phi &0\\
0\ar[r]& D'  \ar[r] & {\C }' \ar[r]^{p'} &{\B }'\ar[r] &0}$$
consists of functors $\phi$ and $\phi _0$, such that $\phi
p=p'\phi _0$, together with a natural transformation $\phi _1:D\to
D'\circ (\phi^{op}\times \phi)$ such that:
$$\phi _0(f_0+\alpha)=\phi _0(f_0)+\phi _1(\alpha)$$
for all $f_0:c\to d$ in $\C $ and $\alpha$ in $D(c,d)$.
\par\medskip
We now list properties of linear extensions relevant to our
problem.
\begin{enumerate}
\item If $\B $ is  a small category, there is \cite[IV.6]{BW} a
canonical bijection
$$   M ({\B }, D) \cong H^2 ({\B }, D). $$
from the set of equivalence classes of linear extensions of $\B$
by $D$ and the second cohomology group $H^2 ({\B },D)$ of $\B $
with coefficients in $D$. \item The functor $p$ reflects
isomorphisms and yields a bijection on the sets of isomorphism
classes ${\rm Iso}({\C })\cong {\rm Iso}({\B })$. \item Let
$(\phi_1,\phi _0,  \phi )$ be a morphism of linear extensions.
Suppose that $\phi _1(c,d)$ is an isomorphism for any $c$ and $d$
in $\B$. Then $\phi $ is an equivalence of categories if and only
if $\phi _0$ is an equivalence of categories. \item If $\B $ is an
additive category and $D$ is a biadditive bifunctor, then the
category $\C $ is additive \cite[Proposition 3.4]{JP}.
\end{enumerate}
\begin{Prop}
Let $D$ be the bifunctor defined on $\B$ by:
$$D(B,B'):=\Hom_{\A}(B/i_*i^!B,i_*i^!B')\ .$$
The category $\A$ is a linear extension of $\B$ by $D$.
\end{Prop}
{\it Proof}. It reduces to the following. Two maps of extensions:
$$\xymatrix{
0\ar[r]& U \ar[r]\ar[d]&\ar[d]_f A\ar[d]^g\ar[r] & X \ar[r]\ar[d]&0\ \ \\
0\ar[r]& U'\ar[r]& A' \ar[r] & X' \ar[r]&0}
$$
agree on the side vertical arrows if and only if their difference
$f-g$ factors through a map in the group $\Hom (X,U')$.\rdg The
results of Section \ref{717} shows that the categories $\A'$,
$\A''$ and the functors $i^*j_*$, $i^!j_!$ of the recollement
situation determine the category $\Ker i^!$. We now show that it
does determine the bifunctor $D$ as well. For an object $B$ in
$\B$, let $((X,V,\alpha),U)$ be its image under the composite:
$$\B\simeq\G\to\Ker i^!\times\A'\simeq\Su(i^*j_*)\times\A'\ .$$
That is: $X=j^*A$, $V=i^*A$, for $A=B/i_*i^!B$, $U=i^!B$. Then:
\begin{equation}\label{22102}
D(B,B'):=\Hom_{\A}(A,i_*U')=\Hom_{\A'}(i^*A,U')=\Hom_{\A'}(V,U')\
.\end{equation}
\section{A comparison theorem}
We have seen in Section \ref{contrex} an example of a comparison
functor which is not an equivalence of categories.
However, a comparison functor $E$ indeed yields an equivalence
from $\Ker (i^*:\A_1\to\A')$ to $\Ker (i^*:\A_2\to \A')$,
and similarly for $\Ker i^!$.
If $E$ is an equivalence of categories, then clearly $E$ commutes
with the derived functors $\R ^{\bullet}i^!$ and $\L _{\bullet}i^*$.
This observation leads to the following definition.
\begin{De}
Let $(\A',\A_1,\A'')$ and $(\A',\A_2,\A'')$ be two recollement
situations.
Assume that the categories $\A_1,\A_2,\A',\A''$ have enough
projective objects. A comparison functor $E$:
$\A_1\to \A_2$ is left admissible if the following diagram
commutes
$$\xymatrix{\A'\ar[d]_{=}& \Ker i^! \ar[l]_{\L _1i^* }\ar[d]^{E}\\
\A'&\Ker i^!\ar[l]^{ \L _1i^*} }
$$
A right admissible comparison functor is defined similarly
by using the functors $\R ^1i^!$ and the categories $\Ker i^*$.
\end{De}
\begin{Th}\label{comparison} Let $E$ be a comparison functor between categories
with enough injectives and projectives. The
following conditions  are equivalent
\begin{enumerate}
    \item $E$ is right admissible
    \item $E$ is left admissible
    \item $E$ is an equivalence of categories.
\end{enumerate}
\end{Th}
{\it Proof}. It is clear that iii) implies both
conditions i) and ii).  We only show that ii) implies iii). A dual argument
shows that i) implies iii). By Section \ref{linext}, the functor $E$
yields a commutative diagram of linear extensions
$$\xymatrix{0\ar[r] &D_1 \ar[r]\ar[d]
 &\A_1 \ar[r]\ar[d]^E &\B_1\ar[r]\ar[d]& 0 \\
0\ar[r] &D_2 \ar[r] &\A_2 \ar[r] &\B_2\ar[r]& 0 }
$$
First we show that $E$ yields an equivalence of categories
$\B_1\to \B_2$. By Section \ref{linext} it suffices to show that
$E$ yields an equivalence $\G_1\to \G_2$. When there are enough
projectives, $E$ yields an equivalence on $\Ker i^!$ (Theorem \ref{qartuli}).
The induced map
$$\Ext^1_{\A_1}(A,i_*U)\to \Ext^1_{\A_2}(E(A),i_*U)$$
is an isomorphism for $U$ in $\A'$ and $A$ in $\Ker i^!$,
thanks to Proposition \ref{spectral} and the five-lemma.
Once $\B_1$ and $\B_2$ are identified, we use the computation
(\ref{22102}) to conclude that the morphism of bifunctors $D_1\to
D_2$ is an isomorphism. The rest is a consequence of the
properties of linear extensions of categories. \rdg
\section{Recollement pr\'e-h\'er\'editaire}
\subsection{{\qh} recollement}
\begin{De}
A recollement situation with enough projectives is {\qh} if
for any projective object $V$ of the category $\A'$:
$$(\L _2i^*)(i_*V)=0\ .$$
\end{De}
\begin{Prop}\label{12324} In a {\qh} recollement situation:
$(\L _2i^*)i_*=0$.
\end{Prop}
{\it Proof}. By (\ref{88}) the functor $(\L _2i^*) i_*$ is right
exact. If it vanishes on projective objects, it vanishes on all
objects. \rdg
\begin{Le}\label{24122}  In a {\qh} recollement situation
there is an isomorphism of functors
$$(\L _1i^*) j_*\cong i^!j_!.$$
\end{Le}
{\it Proof}. Apply the functor $i^*$ to the short exact sequence:
$$0\to j_{!*}\to j_*\to i_*i^*j_*\to 0\ .$$
By (\ref{88}), $\L _1i^*$ vanishes on $i_*i^*j_*$, and by
hypothesis $\L _2i^*$ vanishes on $i_*i^*j_*$. Hence the long
exact sequence for left derived functors yields an isomorphism:
$(\L _1i^*) j_{!*}\cong (\L _1i^*) j_*$. The result follows by
(\ref{1142}). \rdg
\begin{Th}\label{12325} Let $(\A',\A,\A'')$ and $(\A',\B,\A'')$ be two
{\qh} recollement situations and let
 $E:\A\to \B$ be a comparison functor. Then $E$ is admissible and hence is an equivalence of categories.
\end{Th}
{\it Proof}. We have to prove that $\L _1i^*$ has the same value
on $A$ and $EA$, provided that $i^!A=0$.
For such an $A$, there is a short exact sequence (\ref{cokeretaonly}).
Applying the functor $i^*$ results in an exact sequence:
$$\L_2i^*(\Coker\eta_A)\to \L_1i^*(A)\to \L_1i^*(j_*j^*A)\to \L_1i^*(\Coker\eta_A)$$
whose right term cancels by Proposition \ref{etaonKer} and (\ref{88}),
and whose left term cancels by Proposition \ref{12324}.
This gives an isomorphism: $\L _1i^*(A)\cong (\L _1i^*)j_*j^*(A)$.
Lemma \ref{24122} finishes the proof. \rdg
\subsection{MacPherson-Vilonen recollements}
The following proposition is a formalized version of
the construction of projectives in \cite[Proposition 2.5]{MiV}.
\begin{Prop}
Let $\A(F\buildrel \xi \over\rightarrow G)$ be a
Mac-Pherson-Vilonen recollement. Assume further that the left
exact functor $G$ has a left adjoint $G^*$. Then the exact functor
$r$ has a left adjoint $r^*$defined by:
$$r^* V = (G^* V , FG^* V \oplus V ,(1,0),\xi_{G^*V}\oplus\eta_V)$$
where in this formula $\eta$ denotes the unit of adjonction:
$id_{\A'}\to GG^*$. In particular, there is a short exact
sequence:
 \begin{equation}\label{r*}
 0 \to j_! G^* \to r^* \to i_* \to 0\ .
 \end{equation}
\end{Prop}
\textit{Proof.}
 Necessarily, $j^*r^* =(rj_!)^*=G^*$. Then check. \rdg
\begin{Prop}\label{vf}
Every MacPherson-Vilonen recollement with enough projectives is {\qh} .
\end{Prop}
\textit{Proof.} Apply the functor $i^*$ to the short exact sequence
(\ref{r*}). Part of the resulting long exact sequence is an exact
sequence:
$$(\L_2 i^*)r^*\to(\L_2 i^*)i_*\to(\L_1 i^*)j_!G^*\ ,$$
whose right term cancels by (\ref{pirvelicarmoebuli}).
To conclude, if $P$ is a projective in $\A'$, then $r^*P$ is a
projective in $\A$, because $r^*$ is left adjoint to an exact
functor.
\rdg
This leads to the following characterization of MacPherson-Vilonen recollements.
A special case appeared in \cite[Proposition 2.6]{V}
\begin{Th}\label{MV}
A recollement situation of categories with enough projectives
 is isomorphic to a MacPherson-Vilonen
construction if and only if the recollement is {\qh}
and there exists an exact functor $r$: $\A \to \A'$ such that
$r\circ i_*=Id_{\A '}$.
\end{Th}
{\it Proof}. Consider a recollement with such
an exact retraction functor $r$. The natural transformation $N$: $j_!\to j_*$ yields
a transformation $rN$ from the right exact functor $r  j_!$ to the
left exact functor $r  j_*$. Thus we can form the
MacPherson-Vilonen construction $\A(r  j_!\buildrel {rN}
\over\rightarrow r  j_*)$. We define a functor $E$: $\A\to \A(r
j_!\buildrel {rN} \over\rightarrow r  j_*)$ by:
$$E(A)=(j^*(A),r(A), r(\epsilon _A), r(\eta _A)).$$ One
checks with Section \ref{macp} and (\ref{norm}) that $E$ is a comparison functor.
By Proposition \ref{vf}, $\A(rN)$ is {\qh}. If $\A$ is also {\qh}, Theorem
\ref{12325} applies.\rdg
\textsc{Remark.} Similarly one can define
\emph{pre-cohereditary} recollements by the condition $\R ^2
i^!(i_*V)=0$ for any injective $V$ in $\A'$. We leave to the
reader to dualize the above results.
\subsection{The case when $i^*j_*=0$ or $i^!j_!=0$}
In this section, we characterize the recollements
$\A={\A '}\rtimes _F{\A ''}$ of Section \ref{semidirect}.
\begin{Prop}\label{oripiroba}
For a recollement with enough projectives, the following are
equivalent:
 \begin{enumerate}
 \item The functor  $i^*$ is exact.
 \item $i^!j_!=0$.
 \end{enumerate}
 Dually, for a recollement with enough injectives, the following  are equivalent:
\begin{enumerate}
 \item The functor  $i^!$ is exact.
 \item $i^*j_*=0$.
 \end{enumerate}
\end{Prop}
{\it Proof}. We prove the second assertion. Assume that $i^!$
is exact. Applying $i^!$ to the epimorphism $j_*\to i_*i^*j_*$
gets an epimorphism $0=i^!j_*\to i^!i_*i^*j_*\cong i^*j_*$.
\par
Assume conversely that $i^*j_*=0$ and suppose that the recollement
has enough injectives. We first prove that $\R ^1i^!(A)=0$ when
$i^!A=0$. By Proposition \ref{etaonKer}, if $i^!A=0$, there is an
epimorphism $j_*j^*A\to \Coker\eta_A\cong i_*(\R ^1i^!)(A)$.
Applying the right exact functor $i^*$, we get an epimorphism
$i^*j_*j^*(A)\to (\R ^1i^!)(A)$.
\par
Next, we apply $i^!$ to the short exact sequence (\ref{eta}). It
yields an exact sequence:
$$0\to i^!\buildrel {\simeq}\over\rightarrow i^!\to i^!\Im\eta\to
(\R ^1i^!)i_*i^!\to \R ^1i^! \to(\R ^1i^!)\Im\eta\ .$$ By
(\ref{88}), $(\R ^1i^!)i_*i^!=0$, so that: $i^!\Im\eta=0$. It
results that $(\R ^1i^!)\Im\eta=0$ as well, and finally that
$R^1i^!=0$. \rdg
As an application we recover \cite[Proposition 2.4]{A}.
\begin{Prop}\label{daxasiateba} Every recollement situation with enough projectives,
such that: $i^!j_!=0$, is equivalent to $\A'\ltimes_{i^*j_*}\A''$.
Dually, every recollement situation with enough injectives, such
that: $i^*j_*=0$, is equivalent to $\A'\rtimes_{i^!j_!}\A''$.
\end{Prop}
{\it Proof}. When the recollement has enough projectives, Theorem
\ref{MV} applies for $r=i^*$. \rdg
\begin{Co} Let $\A',\A,\A''$ be a recollement situation
with enough projective or enough injectives. If the
norm $N:j_!\to j_*$ is an isomorphism, then $\A\cong \A'\x \A''$.
\end{Co}
{\it Proof}.  By Proposition \ref{birtvi}: $i^*j_*=i^!j_!=0$. Then
we apply Proposition \ref{daxasiateba}. \rdg
\par
\vskip2cm {\sc Acknowledgements}.
The second author would like to thank University of
Nantes for hospitality and support.

\Addresses
\end{document}